\documentclass[12pt]{article}
\usepackage{amsfonts, graphicx, amsmath, verbatim}

        \renewcommand{\Im}{\mathop{\rm Im\,}\nolimits}

        \newcommand{\C}{\mathbb C}

        \newcommand{\dd}{\partial}

        \newtheorem{theorem}{Theorem}

        \newtheorem{lemma}{Lemma}

        \newtheorem{corollary}{Corollary}
        \newcommand{\be}[1]{\begin{equation}\label{#1}}
        \newcommand{\ee}{\end{equation}}
        \newcommand{\myref}[1]{$(\ref{#1})$}

        \newcounter{example}
        
\title{A note on a harmonic measure estimate 
      and a conjecture of J.~Velling}
\author{Alexander Fryntov}
\date{December 5, 2021}

\begin{document}
\maketitle

\begin{abstract}
	Suppose that finitely many disjoint open arcs have been selected
	on the unit circle, each of length less than $\pi$. Let $L_0$
        be a longest among them. One can treat the unit disk as a hyperbolic 
        plane in 
        the Poincare disk model. From this perspective each arc $L$ of the selected set
        determines a hyperbolic half-plane bounded by the 
        geodesic curve $I$ joining endpoints of the arc $L$. Remove from the unit disk 
        all these hyperbolic half-planes. The remaining domain is simply 
        connected and contains the origin. Now map this domain conformally onto the unit
	disk so that the origin stays fixed. After this map, the boundaries of the
	hyperbolic half-planes appear as disjoint arcs on the unit circle.
	Let  $I'_0$ be the conformal image of the boundary of the hyperbolic half 
        plane determined
        by the arc $L_0$.

        In \cite{velling} J.~Velling conjectured that 
           \be{eq:velling_conj}
           |I'_0| \ge |L_0|\,,
           \ee
	where $|\cdot|$ stands for the Euclidean length. He proved some conditional theorems 
        based on validity of this conjecture. 

	In this note we prove a theorem which implies the J.~Velling conjecture
	and thus converts the conditional theorems in \cite{velling}
	into true ones.
\end{abstract}

\subsubsection{Velling and Basic basic model domains.}
All plane objects in this note are assumed embedded in the complex plane $\C$.
So we use standard complex notations $\Delta$ for the unit  
disk and $\dd\Delta$ for its boundary. The set $\nu E$ is defined as 
     $$
     \nu E:= \{\nu z:\, z\in E\}\,,\quad \nu\in \C,\ E\subset \C\,.
     $$
Let $E$ be a domain whose boundary consists of a finite set of smooth Jordan arcs,
and let $K$ be a closed subset of $\dd E$.
Notations  $\omega(z, K, E)$ and $g(z,\xi, E)$ are reserved 
for harmonic measures and Green functions. We remind that the 
harmonic measure of $K$ at $z\in \C$ relative to $E$ is a harmonic
in $z\in E$ function that satisfies the boundary condition
    $$
    \omega(z, K, E) =
    \left\{
    \begin{array}{lll}
    1 & : & z\in K \\
    0 & : & z\in \dd E\setminus K
    \end{array}
    \right.
    $$
If $E\subset E_1$, $K\subset \dd E$, and $K\subset K_1\subset \dd E_1$ then 
the Maximum Principle for harmonic functions implies
   $$
   \omega(z, K, E)\le \omega(z, K_1, E_1)\,,  
   $$ 
and the inequality is strict unless $E=E_1$ and $K=K_1$.
This inequality is known as the Domain Extension Principle.

The Green function of $E$ is
   $$
   g(z,\xi, E) = \log\frac{1}{|z-\xi|} + u_\xi(z)\,,\quad z,\xi\in E\,,
   $$
where $u_\xi(z)$ is harmonic in $z\in E$ with the boundary condition
   $$
   u_\xi(z) = \log|z-\xi|\,,\quad z\in\dd E\,.
   $$
The Maximum Principle and the definition imply that $g$ is positive as $z\in E$ and vanishes
as $z\in \dd E$. 

It is known that
  \be{eq:green_sym}
  g(z,\xi, E) = g(\xi, z,E)\,,\quad z,\xi\in E\,,
  \ee
and the harmonic measure can be evaluated it terms of Green function as follows
  \be{eq:harm_via_green}
  \omega(z,K, E) = \frac{1}{2\pi} \int_{K} \frac{\dd g(z,\xi, E)}{\dd n}\,|d\xi|\,,
  \ee
where $\cdot/\dd n$ is used for inward normal derivative in $\xi$ relative to the domain $E$\,.

We reserve letter $L$ for closed arcs of the unit circle $\dd \Delta$. 
The letter $\Lambda$ is used to denote respective angles
    \be{eq:def_lambda}
    \Lambda = \{z\ne 0:\, z/|z|\in L\}\,.
    \ee
Let $L_k$ be a finite set of arcs of opening less than $\pi$ with disjoint interiors, and 
$\bigcup L_k = \dd \Delta$. Suppose that $L_0$ is a longest among them
Let $\theta_k$ be the center of respective $L_k$\,.
For definiteness we assume that $\theta_0=1$\,. Since the length of every arc $L_k$
does not exceed that of $L_0$, then there exists $\eta_k$  ($|\eta_k|=1$,  $\Im \eta_k\ge 0$)
such that
    \be{eq:def_eta}
    L_k = \theta_k\{\eta_k L_0\cap \bar\eta_k L_0\}\,.
    \ee

Now we describe the {\em Velling model} domain $D$. For every arc $L_k$
take the circular arc {\em orthogonal} to $\dd\Delta$ that joins the endpoints of $L_k$.
Denote this arc by $I_k$\,. Define $C_k$ as the closed set bounded by the pair of arcs
$L_k$ and $I_k$. We call the set
  \be{eq_def_velling_domain}
  D := \Delta\setminus (\cup_k C_k)
  \ee
the {\em Velling model} domain supported by the arcs $L_k$\,.

    \begin{figure}[htb]
    \leavevmode
    \makebox[.01\textwidth]{}
    \centering
    \includegraphics[scale=0.65]{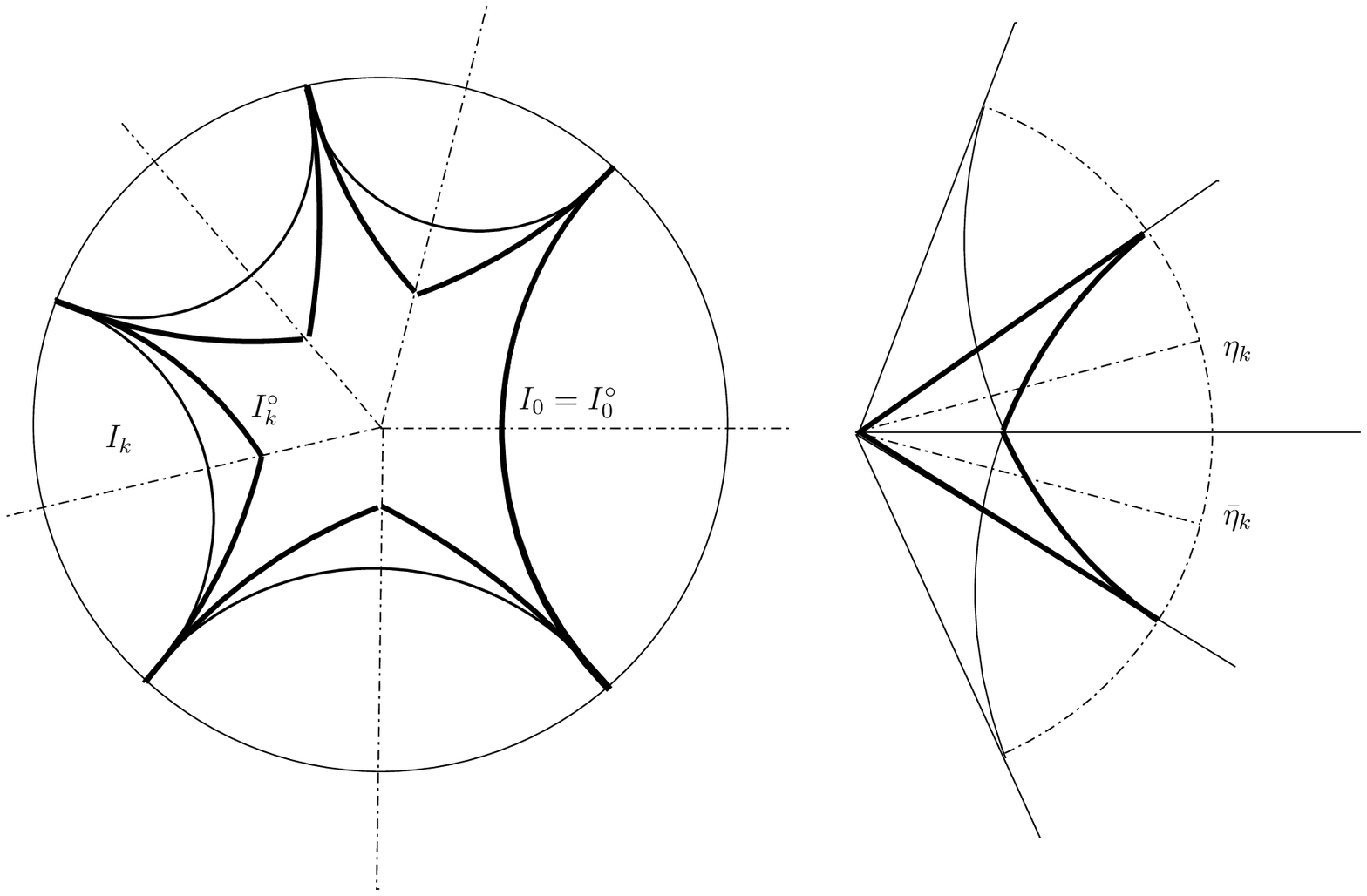}
    \caption{}
    \end{figure}

Since $I_0'$ is the image of the arc $I_0$ under a conformal map of $D$ onto $\Delta$
that leaves the origin intact, then $|I_0'|= 2\pi\omega(0, I_0, D)$. Thus, 
inequality~\myref{eq:velling_conj} becomes
   \be{eq:velling_harm}
   \omega(0,I_0, D)\ge \omega(0, L_0, \Delta)\,.
   \ee
In our construction the arcs $L_k$ {\em fill} 
the whole unit circle, while in the
original Velling conjecture it is not required. The Domain Extension Principle shows
that removal of some arcs (except  $L_0$) from the set $L_k$ cannot reduce the harmonic
measure in LHS of~\myref{eq:velling_harm}. So this extra condition imposed onto arcs
$L_k$ does not restrict generality.


Now we describe the {\em basic} model domain $D^\circ$.  
Let $L_k$ be the same arcs as in the previous model.
A simply connected domain $D^\circ\subset \Delta$ is called {\em basic} model domain
if
  $$
  D^\circ _k = \theta_k(\eta_k D^\circ_0\cap \bar\eta_k D^\circ_0)\,,
  \quad
  D^\circ_k := D^\circ \cap \Lambda_k\,,
  $$ 
where $\theta_k$ (centers) and $\eta_k$ (deviations) are the same as in~\myref{eq:def_eta},
while the arc 
  $$
  I^\circ_0 := \dd D^\circ \cap \Lambda_0 
  $$
satisfies the special properties  
  \begin{itemize}
  \item[a)] $I_0^\circ$ is symmetric relative to the positive ray;
  \item[b)] $\min\{|\xi|:\, \xi\in I_0^\circ\} = \rho_0 < 1$ is attained at $\xi = 1$; 
  \item[c)] for every $\rho\in (\rho_0, 1)$ the intersection $I^\circ_0\cap \{|z|=\rho\}$
            consists exactly of two points.
  \end{itemize}
We call $I^\circ_0$ the {\em basic} arc of the domain $D^\circ$\,.

Notice that if $D$ is the Velling model domain then the arc $I_0$ satisfies 
all the properties to be a {\em basic} arc. The basic domain $D^\circ$ with the basic arc 
$I^\circ_0 := I_0$ is evidently satisfy the property
   $$
   D^\circ \subset D\,.
   $$ 
The Domain Extension Principle then implies that
  \be{eq:comp_d_d_knot}
  \omega(0, I_0, D) \ge \omega(0, I_0, D^\circ)\,.
  \ee

\subsubsection{Harmonic Measure and Main Theorem.}
Velling's conjecture follows from
  \begin{theorem}\label{main} 
      Let $D^\circ$ be a basic model domain then
      $$
      \omega(0, I^\circ_0, D^\circ) \ge \omega(0, L_0, \Delta)\,.
      $$
      Moreover, the inequality is strict unless all the arcs have the same length. 
  \end{theorem} 

Let 
   \be{eq:def_omega_knot} 
   \omega_0 := \omega(0, L_0, \Delta)\,.
   \ee
Consider an auxiliary domain
   \be{eq:def_omega}
   \Omega := \{f(z):\, z\in D^\circ_0\} \cup \{0\}\,,
   \ee   
where $f(z)$, $z\in \C \setminus (-\infty,0]$ is the branch of the function $z^{1/\omega_0}$
such that $f(1)=1$\,.

We need some facts on the Green function $g(z, \xi, \Omega)$ which follow from 
\begin{lemma} Let $\Omega$ be the domain defined above and let $g(z,\xi,\Omega)$
        be its Green function. Then
            \be{eq:g_sym}
            g(z, 0, \Omega) = g(\bar z, 0, \Omega)\,,\quad z\in \Omega\,,
            \ee
        and
            \be{eq:positive_derivative}
            \frac{\dd g(\rho e^{it}, 0, \Omega)}{\dd t} > 0\,,
            \ee
        for every $z=\rho e^{it}\in \Omega_+ = \Omega\cap\{\Im z > 0\}$\,.
\end{lemma}

This lemma can be found in the paper by A.~Baernstein~\cite[p.154, Corollary]{baernstein}. 
 
Now we give a sketch of this proof. Consider the function 
$u(z)$ in $\Omega_+ = \Omega\cap\{\Im z> 0\}$ defined as follows
     $$
     u(\rho e^{it}) = \frac{\dd g(\rho e^{it})}{\dd t}\,.
     $$ 
Simple evaluation of Laplacians in polar coordinates implies that $u(z)$ is
harmonic. Property a) of the basic arc $I^\circ_0$ implies equation~\myref{eq:g_sym}.
Properties b) and c) imply that
boundary values of $u$ are positive on 
$\dd \Omega \cap \{\Im z > 0\}$ and zero on the rest of $\dd \Omega_+$.
The Maximum Principle then implies that $u\ge 0$. Since $u$ is not identical zero
then it is strictly positive.

\subsubsection{Proof of the Theorem.}
Let $\phi$ be a function such that 
  $$
  u(z):= \phi(z) -g(z,0,D^\circ)
  $$ 
is continuous on
the closure of $D^\circ$ and vanishes on the boundary $\dd D^\circ$. 
Suppose that $u(z)$ is subharmonic
on $D^\circ$ and 
    $$
    \int_{I_0^\circ} \frac{\dd \phi(\xi)}{\dd n}\,|d\xi| = 2\pi\omega_0\,.
    $$
To prove the theorem it is enough to show that such a function $\phi$ does exist.

Indeed, by the Maximum Principle $u(z)\le 0$ on $D^\circ$ and the inequality is 
strict everywhere on $D^\circ$ unless $u(z)\equiv 0$. 
It implies that $\dd u(\xi)/\dd n \ne  0$ at every regular point of  
$\dd D^\circ$. Since 
   $$
   \int_{I^\circ_0}\frac{\dd g(\xi, 0, D^\circ)}{\dd n}\,|d\xi| = 
      2\pi \omega(0, I^\circ_0, D^\circ)\,,
   \quad
   \int_{I^\circ_0}\frac{\dd \phi(\xi)}{\dd n}\,|d\xi| = 2\pi \omega(0, L_0,\Delta)\,,
   $$
then
   $$
   \omega(0, I^\circ_0, D^\circ) \ge \omega(0, L_0,\Delta)\,,  
   $$
and the inequality is strict unless 
   \be{eq:g_phi_equal}
   \phi(z) \equiv g(z,0,D^\circ)\,. 
   \ee

Let us construct such a function $\phi(z)$ on $D^\circ$.
In the sector $D^\circ_0$ define the function
   \be{eq:def_phi_knot}
   \phi_0(z) : = \omega_0\,g(f(z), 0, \Omega)\,,\quad z\in D^\circ_0\,.
   \ee
where $f(z)$ is the branch of $z^{1/\omega_0}$ such that $f(1)=1$\,.
Notice that
  \be{eq:on_arc}
  \int_{I^\circ_0}\frac{\dd \phi_0(\xi)}{\dd n}\,|d\xi| = 2\pi\omega_0\,.
  \ee
Take any $\eta\in L_0\cap \{\Im z > 0\}$ and define 
   \be{eq:def_psi}
   \psi_\eta(z) := \phi_0(\eta z)\lor\phi_0(\bar\eta z)\,\quad 
   z\in \eta D^\circ_0\cap \bar \eta D^\circ_0\,.  
   \ee
According to the lemma $\psi_\eta(z) = \psi_\eta(\bar z)$, 
and  the maximum is $\phi_0(z\eta)$ as $\Im z \ge 0$ 
and  $\phi_0(z\bar\eta)$ as $\Im z \le 0$. Therefore,  
  \be{eq:def_by_parts}
  \psi_\eta(z) = \left\{
     \begin{array}{lll}
       \phi_0 (\eta z) &:& \Im z\ge 0\\
       \phi_0 (\bar\eta z)     &:& \Im z < 0\,.
     \end{array}
     \right. \,.
  \ee
Using $\theta_k$ and $\eta_k$ (being the center and deviation of $L_k$)
define
   $$
   \phi_k(z) := \psi_{\eta}(\bar\theta z)\,,\quad
   \text{where}\ \theta=\theta_k,\ \eta=\eta_k,\  z\in D^\circ_k\,. 
   $$
Notice that $\phi_k$ is subharmonic inside of the sector $D^\circ_k$
as the maximum of two harmonic functions.

Notice that if the sets $D^\circ_l$ and $D^\circ_m$ do intersect then their 
intersection is an interval $\gamma_{lm}$, then $\phi_l$ and $\phi_m$ coincide on 
$\gamma_{lm}$. Therefore, there exists a function $\phi$ defined on 
$D^\circ\setminus\{0\}$ and subharmonic inside of every $D^\circ_k$.

Now we will show that the function $\phi$ is {\em harmonic} on $\gamma_{lm}$.
Let $\xi\in\gamma_{lm}$. Select $\delta>0$ so that the ball
$B_\delta(\xi):=\{z:\,|z-\xi|<\delta\}$ contains neither origin nor point 
of the rays $\{z:\,z/|z| = \theta_l\}$ and $\{z:\,z/|z| = \theta_m$.
It follows from~\myref{eq:def_by_parts} and the definitions of $\phi_0$ and $\phi_k$ that
for every $z\in D^\circ_l$
  \be{eq:phi_l_m}
  \phi_l(z) = \phi_m(z^*)\,,
  \ee
where $z^*$ is the reflection of $z$ in $\gamma_{lm}$\,.
On the other hand, the function $\phi_l$ admits a harmonic extension to the whole
ball $B_\delta(\xi)$ such that
  \be{eq:phi_ext}
  \phi_l(z*) = \phi_l(z)
  \ee
(this follows from the definitions of $\phi_0$ and $\phi_l$.)
The equations \myref{eq:phi_l_m} and \myref{eq:phi_ext} then imply that $\phi_m$ 
coincides with the harmonic extension of $\phi_l$ onto $B_\delta(\xi)\cap D^\circ_m$, and hence,
$\phi(z)$ is harmonic at $\xi$\,. 

Thus, the function $\phi$ is subharmonic on $D^\circ\setminus\{0\}$.
The function $u(z) = \phi(z) - g(z,0, D^\circ)$ as it follows
from~\myref{eq:def_phi_knot} is subharmonic on $D^\circ\setminus \{0\}$
and behaves near the origin as 
   $$
   u(z) = o(|\log|z||)\,,\quad z\to 0\,.
   $$
Thus, $u(z)$ is subharmonic at $0$ as well. This completes the proof.


\subsubsection{Remark.}

The following harmonic measure estimate is another
elementary corollary of this theorem (although I have no slightest
idea whether this result is new.)

\begin{corollary} Let $ P_n$ be a bounded domain enclosed by a Euclidean  polygon 
inscribed into the unit circle.
Let $I_0$ be a side of the polygon of largest length and $L_0$ be the
arc supporting this side. If $0\in P_n$ then 
	$$
	\omega(0,I_0, P_n)\ge \omega(0, L_0,\Delta)\,.
	$$
\end{corollary}
\medskip\par
The results of this note was proven when I was a visiting professor at Virginia Tech.
I would like to express my thanks for support and hospitality
of Mathematics Department where the main results was made.

My special thanks are to A.~Solynin and A.~Eremenko
who encouraged me to publish this material.



\begin{thebibliography}{1}

\bibitem{velling}
Velling~John~A.,
\newblock {\it Harmonic measure, infinite kernels, and symmetrization}.
\newblock(English. English summary)
Proc. Amer. Math. Soc. 124 (1996), no. 12, 3739--3743.

\bibitem{baernstein} 
Baernstein~Albert~II.,
\newblock{\it Integral means, univalent functions and
circular symmetrization,}
\newblock{Acta mathematica. {\bf 133}(1974), pp.139--169.}

\end{thebibliography}
\end{document}